\documentclass[a4paper]{amsart}
\usepackage{amssymb} 
\usepackage[T1]{fontenc}
\usepackage[latin1]{inputenc}
\usepackage{amsfonts}
\usepackage{amsxtra}
\usepackage{ae}
\usepackage[all]{xy}
\usepackage{enumerate}
\usepackage{url}
\usepackage{verbatim}
\usepackage{pgf,tikz}
\usepackage{tikz-cd}
\usepackage{stmaryrd}

\usetikzlibrary{arrows,patterns}

\include{diagram}

\newcommand*{\ket}{\rangle}
\newcommand*{\bra}{\langle}
\newcommand*{\ad}{\mathsf{ad}}

\newcommand*{\C}{\mathcal{C}}

\newcommand*{\E}{\mathcal{E}}
\newcommand*{\F}{\mathcal{F}}

\renewcommand*{\H}{\mathcal{H}}

\newcommand*{\I}{\mathcal{I}}
\newcommand*{\J}{\mathcal{J}}

\renewcommand*{\L}{\mathcal{L}}
\newcommand*{\M}{\mathcal{M}}
\newcommand*{\N}{\mathcal{N}}

\renewcommand*{\P}{\mathcal{P}}

\newcommand*{\R}{\mathcal{R}}

\newcommand*{\T}{\mathcal{T}}

\newcommand*{\cop}{\mathsf{cop}}
\renewcommand*{\top}{\mathsf{top}}

\newcommand*{\hyphenAlg}{$-$\mathsf{Alg}}

\newcommand*{\HH}{\mathbb{H}}
\newcommand*{\KH}{\mathbb{K}}
\newcommand*{\LH}{\mathbb{L}}

\newcommand*{\red}{\mathsf{r}}

\newcommand*{\CI}{\mathcal{CI}}
\newcommand*{\CC}{\mathcal{CC}}

\renewcommand*{\max}{\mathsf{f}}

\newcommand{\weights}{\mathbf{P}}

\DeclareMathOperator{\res}{res}
\DeclareMathOperator{\ind}{ind}

\DeclareMathOperator{\id}{id}
\DeclareMathOperator{\ev}{ev}

\numberwithin{equation}{section}
\theoremstyle{change}
\newtheorem{theorem}{Theorem}[section]
\newtheorem{prop}[theorem]{Proposition}
\newtheorem{lemma}[theorem]{Lemma}

\newtheorem{definition}[theorem]{Definition}

\begin{document}

\title{On the assembly map for complex semisimple quantum groups}

\author{Christian Voigt}
\address{School of Mathematics and Statistics \\
         University of Glasgow \\
         University Place \\
         Glasgow G12 8SQ \\
         United Kingdom 
}
\email{christian.voigt@glasgow.ac.uk}


\maketitle

\begin{abstract}
We show that complex semisimple quantum groups, that is, Drinfeld doubles of $ q $-deformations of compact semisimple Lie groups, 
satisfy a categorical version of the Baum-Connes conjecture with trivial coefficients. This approach, based on homological algebra in 
triangulated categories, is compatible with the previously studied deformation picture of the assembly map, and 
allows us to define an assembly map with arbitrary coefficients for these quantum groups. 
\end{abstract} 

\section{Introduction}

If $ G $ is a second countable locally compact group and $ A $ a separable $ G $-$ C^* $-algebra, then the Baum-Connes conjecture for $ G $ with coefficients 
in $ A $ asserts that the assembly map 
$$ 
\mu_A: K^{\top}_*(G;A) \rightarrow K_*(G \ltimes_\red A)
$$
is an isomorphism \cite{BC}, \cite{BCH}. Here $ K^{\top}_*(G;A) $ denotes the topological $ K $-theory of $ G $ with coefficients 
in $ A $, and $ K_*(G \ltimes_\red A) $ is the $ K $-theory of the reduced crossed product.  

For $ G $ a connected Lie group and $ A = \mathbb{C} $ the trivial $ G $-$ C^* $-algebra the assembly map is known to be an 
isomorphism \cite{Lafforguebc}, \cite{CENconneskasparov}, but it is an open problem whether this is true for arbitrary coefficients $ A $ 
in the case of semisimple groups. 

Meyer and Nest have reformulated the Baum-Connes conjecture using the language of triangulated categories and derived 
functors \cite{MNtriangulated}. In their approach, the left hand side of the assembly map is identified with the 
localisation $ \LH F $ of the functor $ F(A) = K_*(G \ltimes_\red A) $ on the equivariant Kasparov category $ KK^G $. 
One of the advantages of the Meyer-Nest picture is that it is very flexible, and applicable in situations where a 
geometric definition of the left hand side of the Baum-Connes assembly map is no longer available.
This has previously been taken advantage of in the setting of discrete quantum groups, leading to explicit $ K $-theory 
computations, see for instance \cite{VVfreeu}, \cite{Voigtqaut}, \cite{FreslonMartoswreath}. 

The present paper is devoted to the construction of an assembly map for complex semisimple quantum groups. 
These locally compact quantum groups, obtained as Drinfeld doubles of $ q $-deformations of compact semisimple Lie groups and 
studied first by Podle\'s and Woronowicz \cite{PWlorentz}, feature in various contexts. 
This includes, among other things, combinatorial quantization of complex Chern-Simons theory \cite{BRLorentz}, and the construction of subfactors 
with property (T), see \cite{Aranospherical}.  

We adopt the framework of Meyer and Nest to the situation at hand. If $ G_q $ is a complex semisimple quantum group 
then it contains the $ q $-deformation $ K_q $ of the corresponding compact semisimple group $ K $ naturally as a quantum subgroup. 
We consider the associated induction and restriction functors and define compactly induced and compactly contractible actions in 
analogy with the case of classical semisimple groups. The analysis in the quantum case is in fact somewhat simpler since $ K_q \subset G_q $ 
is not only compact but also an open quantum subgroup. The localising subcategory of the equivariant Kasparov category $ KK^{G_q} $ generated 
by compactly induced actions is complementary to the category of compactly contractible actions, so that general machinery from homological algebra 
allows us to define an assembly map for arbitrary $ G_q $-$ C^* $-algebras. 
In fact, for technical reasons we will mostly work at the level of continuous fields, assembling the quantum groups $ G_q $ 
for varying deformation parameters into a global object. The corresponding variant of equivariant $ KK $-theory will be developed along the way as far 
as needed for our purposes.

In order to analyse the assembly map we consider a specific projective resolution of the trivial action in $ KK^{G_q} $, 
built out of the Koszul complex for the module of differential forms over the representation ring of $ K_q $.  
According to general results in \cite{Meyerhomalg2}, this leads to a cellular approximation of the trivial action 
in $ KK^{G_q} $. Since $ G_q $ has a unitary braiding, one can also obtain cellular approximations 
for arbitrary actions in this way, using braided tensor products in the sense of \cite{NVpoincare}. 

Our main result is that the assembly map for $ G_q $ with trivial coefficients is an isomorphism. 
This is proved by invoking the results from \cite{MVcqgdeformation} on the deformation of a complex quantum group into its
associated quantum motion group. In this way we provide the first examples of genuine locally compact quantum groups satisfying the 
categorical Baum-Connes property beyond the discrete case. Our argument shows at the same time that the categorical assembly map is 
naturally compatible with the deformation picture developed in \cite{MVcqgdeformation}. 

Let us briefly explain how the paper is organised. Section \ref{secprelim} contains some preliminaries on complex semisimple quantum groups 
and other background material. In section \ref{seccflcqg} we discuss the notion of 
continuous fields of locally compact quantum groups. We explain in particular how Drinfeld doubles of compact 
semisimple Lie groups assemble naturally into continuous fields of locally compact quantum groups by varying the deformation parameter. 
Section \ref{seckkg} develops the basics of equivariant $ KK $-theory for continuous fields of locally compact quantum groups. 
This generalises the corresponding constructions for locally compact quantum groups, corresponding to 
the case that the base space is a point, and groupoid equivariant $ KK $-theory in the sense of \cite{LeGallkkg}
in the case that all fibers are classical groups. We define induction functors for continuous fields of open quantum subgroups, 
and verify that induction is left adjoint to restriction on the level of equivariant $ KK $-theory. 
In section \ref{secbc} we set up the homological machinery for equivariant $ KK $-theory of complex quantum groups in the 
spirit of Meyer and Nest and construct the categorical Baum-Connes assembly map. Moreover we show that 
the assembly map with trivial coefficients is an isomorphism. In order to do so we consider a specific model of the Dirac morphism 
for complex quantum groups, constructed from the Koszul resolution of the representation ring of the classical maximal compact subgroup. 

Let us conclude with some remarks on notation. We write $ \LH(\E) $ for the algebra of adjointable operators on a Hilbert module $ \E $. 
Moreover $ \KH(\E) $ denotes the algebra of compact operators. The closed linear span of a subset $ X $ of a Banach space is denoted 
by $ [X] $. Depending on the context, the symbol $ \otimes $ denotes either the tensor product of Hilbert spaces, the spatial tensor 
product of $ C^* $-algebras, or the spatial tensor product of von Neumann algebras. 
We use standard leg numbering notation for operators defined on multiple tensor products. 

I would like to thank the anonymous referees for their careful reading of the manuscript.

\section{Preliminaries} \label{secprelim}

In this section we collect some background material on locally compact quantum groups in general and complex semisimple quantum groups in particular. 
For more details and further information we refer to \cite{VYcqg}, \cite{MVcqgdeformation}. 

Let $ \phi $ be a normal, semifinite and faithful weight on a von Neumann algebra $ M $. We use the standard notation
$$
\M^+_\phi = \{ x \in M_+| \phi(x) < \infty \}, \qquad \N_\phi = \{ x \in M| \phi(x^* x) < \infty \}
$$
and write $ M_*^+ $ for the space of positive normal linear functionals on $ M $. Assume that
$ \Delta: M \rightarrow M \otimes M $ is a normal unital $ * $-homomorphism. The weight $ \phi $ is called left invariant
with respect to $ \Delta $ if
$$
\phi((\omega \otimes \id)\Delta(x)) = \phi(x) \omega(1)
$$
for all $ x \in \M_\phi^+ $ and $ \omega \in M_*^+ $. Similarly one defines the notion of a right invariant
weight.

\begin{definition} \label{defqg}
A locally compact quantum group $ G $ is given by a von Neumann algebra $ L^\infty(G) $ together with a normal unital $ * $-homomorphism
$ \Delta: L^\infty(G) \rightarrow L^\infty(G) \otimes L^\infty(G) $ satisfying the coassociativity relation
$$
(\Delta \otimes \id)\Delta = (\id \otimes \Delta)\Delta
$$
and normal semifinite faithful weights $ \phi $ and $ \psi $ on $ L^\infty(G) $ which are left and right invariant, respectively.
\end{definition}

Our notation for locally compact quantum groups is intended to make clear how ordinary locally compact groups can be viewed as quantum groups, but 
for a general locally compact quantum group $ G $ the notation $ L^\infty(G) $ is purely formal.
Similar remarks apply to the $ C^* $-algebras $ C^*_\max(G), C^*_\red(G) $ and $ C^\max_0(G), C^\red_0(G) $ associated to $ G $ that we discuss below.
It is convenient to view all of them as different appearances of the quantum group $ G $. 

Let $ G $ be a locally compact quantum group and let $ \Lambda: \N_\phi \rightarrow L^2(G) $ be a GNS-construction for the weight $ \phi $.
We will only consider quantum groups for which $ L^2(G) $ is a separable Hilbert space. One obtains a unitary
$ W_G = W $ on $ L^2(G) \otimes L^2(G) $ by
$$
W^*(\Lambda(x) \otimes \Lambda(y)) = (\Lambda \otimes \Lambda)(\Delta(y)(x \otimes 1))
$$
for all $ x, y \in \N_\phi $. This unitary is multiplicative, which means that $ W $ satisfies the pentagonal equation
$$
W_{12} W_{13} W_{23} = W_{23} W_{12}.
$$
From $ W $ one can recover the von Neumann algebra $ L^\infty(G) $ as the strong closure of the algebra
$ (\id \otimes \LH(L^2(G))_*)(W) $ where $ \LH(L^2(G))_* $ denotes the space of normal linear functionals on $ \LH(L^2(G)) $. Moreover
one has
$$
\Delta(x) = W^*(1 \otimes x) W
$$
for all $ x \in M $. 

The group-von Neumann algebra $ \L(G) $ of the quantum group $ G $ is the strong
closure of the algebra $ (\LH(L^2(G))_* \otimes \id)(W) $ with the comultiplication $ \hat{\Delta}: \L(G) \rightarrow \L(G) \otimes \L(G) $
given by
$$
\hat{\Delta}(y) = \hat{W}^* (1 \otimes y) \hat{W},
$$
where $ \hat{W} = \Sigma W^* \Sigma $ and $ \Sigma \in \LH(L^2(G) \otimes L^2(G)) $ is the flip map. It defines
a locally compact quantum group $ \hat{G} $ which is called the dual of $ G $. The left invariant weight
$ \hat{\phi} $ for the dual quantum group has a GNS-construction $ \hat{\Lambda}: \N_{\hat{\phi}} \rightarrow L^2(G) $,
and according to our conventions we have $ \L(G) = L^\infty(\hat{G}) $. 
We write $ \check{G} $ for the locally compact quantum group obtained from $ \hat{G} $ by equipping the underlying 
von Neumann algebra $ \L(G) $ with the opposite comultiplication 
$ \hat{\Delta}^{\cop} = \sigma \hat{\Delta} $, where $ \sigma $ is the flip map. 

The modular conjugations of the weights $ \phi $ and $ \hat{\phi} $ are denoted by $ J $ and $ \hat{J} $, respectively.
These operators implement the unitary antipodes in the sense that
$$
R(x) = \hat{J} x^* \hat{J}, \qquad \hat{R}(y) = J y^* J
$$
for $ x \in L^\infty(G) $ and $ y \in \L(G) $. Note that $ L^\infty(G)' = JL^\infty(G) J $ and
$ \L(G)' = \hat{J} \L(G) \hat{J} $ for the commutants of $ L^\infty(G) $ and $ \L(G) $. Using $ J $ and $ \hat{J} $ one obtains multiplicative
unitaries
$$
V = (\hat{J} \otimes \hat{J})\hat{W}(\hat{J} \otimes \hat{J}), \qquad \hat{V} = (J \otimes J) W (J \otimes J)
$$
which satisfy $ V \in \L(G)' \otimes L^\infty(G) $ and $ \hat{V} \in L^\infty(G)' \otimes \L(G) $, respectively. 

We will mainly work with the $ C^* $-algebras associated to the locally compact quantum group $ G $. The
algebra $ [(\id \otimes \LH(L^2(G))_*)(W)] $ is a strongly dense $ C^* $-subalgebra of $ L^\infty(G) $
which we denote by $ C^\red_0(G) $. Dually, the algebra
$ [(\LH(L^2(G))_* \otimes \id)(W)] $ is a strongly dense $ C^* $-subalgebra of $ \L(G) $
which we denote by $ C^*_\red(G) $.
These algebras are the reduced algebra of continuous functions vanishing at infinity
on $ G $ and the reduced group $ C^* $-algebra of $ G $, respectively. One has
$ W \in M(C^\red_0(G) \otimes C^*_\red(G)) $, and  restriction of the comultiplications on the von Neumann level
turns $ C^\red_0(G) $ and $ C^*_\red(G) $ into Hopf $ C^* $-algebras in the sense of \cite{NVpoincare}. 

For every locally compact quantum group $ G $ there exists a universal dual $ C^*_\max(G) $ of $ C_0^\red(G) $ and
a universal dual $ C^\max_0(G) $ of $ C^*_\red(G) $, respectively \cite{Kustermansuniversal}.
We call $ C^*_\max(G) $ the maximal group $ C^* $-algebra of $ G $ and
$ C_0^\max(G) $ the maximal algebra of continuous functions on $ G $ vanishing at infinity.
Since $ L^2(G) $ is assumed to be separable the $ C^* $-algebras $ C^\max_0(G), C^\red_0(G) $ and
$ C^*_\max(G), C^*_\red(G) $ are separable. The quantum group $ G $ is called compact if $ C^\max_0(G) $ is unital, and it is called
discrete if $ C^*_\max(G) $ is unital. In the compact case we also write $ C^\max(G) $ and $ C^\red(G) $ instead of $ C^\max_0(G) $ and 
$ C^\red_0(G) $, respectively. 

In general, we have a surjective morphism $ \hat{\pi}: C^*_\max(G) \rightarrow C^*_\red(G) $ of Hopf $ C^* $-algebras
associated to the left regular corepresentation $ W \in M(C_0(G) \otimes C^*_\red(G)) $. Similarly, there is
a surjective morphism $ \pi: C^\max_0(G) \rightarrow C^\red_0(G) $.
We will call the quantum group $ G $ amenable if $ \hat{\pi}: C^*_\max(G) \rightarrow C^*_\red(G) $
is an isomorphism and coamenable if $ \pi: C^\max_0(G) \rightarrow C^\red_0(G) $
is an isomorphism. If $ G $ is amenable or coamenable, we simply write $ C^*(G) $ and $ C_0(G) $
for the corresponding $ C^* $-algebras, respectively. 

An open quantum subgroup $ H $ of a locally compact quantum group $ G $ is given by a surjective normal 
unital $ * $-homomorphism $ \pi: L^\infty(G) \rightarrow L^\infty(H) $ which is compatible with comultiplications \cite{KKSopensubgroups}. 
The corresponding morphism $ \pi: C^\red_0(G) \rightarrow M(C^\red_0(H)) $ takes values in $ C^\red_0(H) $,  
and there exists a central projection $ 1_H \in M(C^\red_0(G)) $ 
such that the map $ \pi $ identifies with $ \pi(f) = 1_H f $. 
Using the projection $ 1_H $ we obtain 
a canonical embedding $ \iota: C^\red_0(H) \rightarrow C^\red_0(G) $ from the identification $ C^\red_0(H) \cong 1_H C^\red_0(G) $. 
This map can be interpreted as extending functions on $ H $ by zero to all of $ G $. 
Moreover we have a canonical inclusion map $ C^*_\red(H) \rightarrow C^*_\red(G) $, compatible with the comultiplications. 

For the definition of actions of locally compact quantum groups on $ C^* $-algebras and their crossed products we refer to \cite{NVpoincare}, 
or to the discussion in section \ref{seckkg} below. 
We shall gather here some results from \cite{Vaesimprimitivity} concerning induced $ C^* $-algebras. 
Assume that $ G $ is a regular locally compact quantum group and that $ H \subset G $ is a closed quantum subgroup. 
In fact, in the sequel we will only be interested in the case that $ H \subset G $ is open. 
From the quantum subgroup $ H \subset G $ one obtains a right coaction $ L^\infty(G) \rightarrow L^\infty(G) \otimes L^\infty(H) $ by right translations 
on the level of von Neumann algebras. The von Neumann algebraic homogeneous space $ L^\infty(G/H) \subset L^\infty(G) $ is defined as the subalgebra of invariants 
under this coaction. If $ \hat{\pi}': \L(H)' \rightarrow \L(G)' $ is the homomorphism $ \hat{\pi}'(x) = \hat{J}_G \hat{\pi}(\hat{J}_H x \hat{J}_H) \hat{J}_G $
induced by $ \hat{\pi}: \L(H) \rightarrow \L(G) $, then
$$
I = \{v \in \LH(L^2(H), L^2(G)) \mid v x = \hat{\pi}'(x) v \;\text{for all}\; x \in \L(H)' \}
$$
defines a von Neumann algebraic imprimitivity bimodule between the von Neumann algebraic crossed product $ G \ltimes L^\infty(G/H) $ and $ \L(H) $.
There is a $ C^* $-algebraic homogeneous space $ C^\red_0(G/H) \subset L^\infty(G/H) $ and
a $ C^* $-algebraic imprimitivity bimodule $ \I \subset I $ which implements a Morita equivalence
between $ G \ltimes_\red C^\red_0(G/H) $ and $ C^*_\red(H) $. Explicitly, we have
$$
\I \otimes L^2(G) = [\hat{V}_G (I \otimes 1)(\id \otimes \hat{\pi})(\hat{V}_H^*)(C^*_\red(H) \otimes L^2(G))].
$$
If $ A $ is an $ H $-$ C^* $-algebra with coaction $ \alpha: A \rightarrow M(C_0^\red(H) \otimes A) $ then the induced $ G $-$ C^* $-algebra $ \ind_H^G(A) $ 
is defined in \cite{Vaesimprimitivity} by first constructing $ G \ltimes_\red \ind_H^G(A) $, and then applying a quantum version of Landstad's theorem 
characterising reduced crossed products. The crossed product of the induced $ C^* $-algebra $ \ind_H^G(A) $ is 
$$
G \ltimes_\red \ind_H^G(A) = [(\I \otimes 1) \alpha(A) (\I^* \otimes 1)], 
$$  
and the Hilbert module implementing the Morita equivalence $ G \ltimes_\red \ind_H^G(A) \sim_M H \ltimes_\red A $ 
can be concretely described as 
$$ 
\J = [(\I \otimes 1)\alpha(A)] \subset M(\KH(L^2(H), L^2(G)) \otimes A). 
$$
We note that \cite{Vaesimprimitivity} supposes strong regularity of the quantum group $ G $, but for the above constructions 
on the reduced level regularity is sufficient. 

Let us now assume in addition that $ H \subset G $ is an open quantum subgroup. 
Since $ 1_H \in M(C_0(G/H)) $ we have $ 1_H \I \subset \I $,  
which implies $ p = 1_H \otimes \id \in M(\J \J^*) = M(G \ltimes_\red \ind_H^G(A)) $. 
Using that $ 1_H \I $ can be identified with $ C^*_\red(H) $, viewed as subspace of $ \LH(L^2(H), L^2(G)) $ via the canonical 
inclusion $ i: L^2(H) \rightarrow L^2(G) $, 
one checks $ p (G \ltimes_\red \ind_H^G(A)) p = H \ltimes_\red A $ 
and $ (G \ltimes_\red \ind_H^G(A)) p (G \ltimes_\red \ind_H^G(A)) = (G \ltimes_\red \ind_H^G(A)) $. 
In other words, the projection $ p $ exhibits $ H \ltimes_\red A $ as a full corner of $ G \ltimes_\red \ind_H^G(A) $ in this case. 

After these general considerations let us now briefly review the definition of complex semisimple quantum groups, referring to \cite{VYcqg} for the details. 
We start with a simply connected semisimple complex Lie group $ G $, its associated Lie algebra $ \mathfrak{g} $, and a positive real deformation 
parameter $ q \neq 1 $. Fix a Cartan subalgebra $ \mathfrak{h} \subset \mathfrak{g} $, and 
let $ \mathfrak{k} \subset \mathfrak{g} $ be the Lie algebra of the maximal compact subgroup $ K \subset G $ with maximal torus $ T $ 
such that $ \mathfrak{t} \subset \mathfrak{h} $, where $ \mathfrak{t} $ is the Lie algebra of $ T $.  
We fix simple roots $ \alpha_1, \dots, \alpha_N $, their corresponding coroots, and denote by $ \weights $ the set of weights for $ \mathfrak{g} $. 
The set $ \weights^+ \subset \weights $ of dominant integral weights is the set of all non-negative integer combinations of the fundamental 
weights $ \varpi_1, \dots, \varpi_N $.  

The $ C^* $-algebra $ C(K) $ of continuous functions on the maximal compact subgroup $ K $ of $ G $ can be described 
as completion of the algebra of matrix coefficients of all finite dimensional representations of $ \mathfrak{g} $. 
In a similar way one constructs the $ C^* $-algebra $ C(K_q) $ as completion of the algebra 
of matrix coefficients of all finite dimensional weight modules for the quantized universal enveloping algebra $ U_q(\mathfrak{g}) $ 
associated with $ \mathfrak{g} $. 

The irreducible finite dimensional weight modules of $ U_q(\mathfrak{g}) $ are parametrised by their highest weights in $ \weights^+ $ as 
in the classical theory. We will write $ V(\mu)^q $, or simply $ V(\mu) $, for the module associated to $ \mu \in \weights^+ $. 
We fix the $ * $-structure on $ U_q(\mathfrak{g}) $ as in \cite{VYcqg} and note that the representations $ V(\mu)^q $ for $ \mu \in \weights^+ $ are unitarizable 
with respect to this $ * $-structure. 
In fact, one can identify the underlying Hilbert spaces $ V(\mu)^q $ with the Hilbert space $ V(\mu)^1 = V(\mu) $ of the corresponding irreducible representation of 
$ U(\mathfrak{g}) $ in a natural way, compare \cite{NTKhomologydirac}. 

By construction, the $ C^* $-algebra $ C(K_q) $ contains the canonical dense Hopf $ * $-algebra of matrix coefficients, 
which has a linear basis given by the elements $ u^\mu_{ij} = \bra e^\mu_i|\bullet| e^\mu_j \ket $, 
where $ \mu \in \weights^+ $ and $ e^\mu_1, \dots, e^\mu_n $ is an orthonormal basis of $ V(\mu) $ consisting of weight vectors. 
In terms of matrix elements, the Haar state $ \phi: C(K_q) \rightarrow \mathbb{C} $ is determined by 
$$
\phi(u^\mu_{ij}) = 
\begin{cases}
1 & \text{if } \mu = 0 \\
0 & \text{else. }
\end{cases}
$$
We will write $ \Delta = \Delta^q: C(K_q) \rightarrow C(K_q) \rightarrow C(K_q) $ for the comultiplication of $ C(K_q) $. 
Dually, the group $ C^* $-algebra $ C^*(K_q) $ is the $ c_0 $-direct sum 
$$
C^*(K_q) = \bigoplus_{\mu \in \weights^+} \KH(V(\mu)),   
$$
and we write $ \hat{\Delta} = \hat{\Delta}^q: C^*(K_q) \rightarrow M(C^*(K_q) \otimes C^*(K_q)) $ for its comultiplication. 
The algebra $ C^*(K_q) $ contains the algebraic direct sum of the matrix algebras $ \KH(V(\mu)) $ 
as a dense subalgebra, and one can choose a basis $ \omega^\mu_{ij} $ for each $ \KH(V(\mu)) $ dual to the matrix coefficients $ u^\mu_{ij} $.  
With these choices, the fundamental multiplicative unitary $ W \in M(C(K_q) \otimes C^*(K_q)) $ can be written as the strong$^*$ convergent sum 
$$
W = \sum_{\mu \in \weights^+} \sum_{i,j = 1}^{\dim(V(\mu))} u^\mu_{ij} \otimes \omega^\mu_{ij},
$$
see \cite{VYcqg}. 

The complex semisimple quantum group $ G_q $ is the Drinfeld double $ G_q = K_q \bowtie \hat{K}_q $, with underlying $ C^* $-algebra
$$
C_0(G_q) = C(K_q) \otimes C^*(K_q)
$$
and comultiplication 
$$
\Delta_{G_q} = (\id \otimes \sigma \otimes \id)(\id \otimes \ad(W) \otimes \id)(\Delta \otimes \hat{\Delta}).
$$
Here $ W $ denotes the multiplicative unitary from above. 

By definition, a unitary representation of $ G_q $ on a Hilbert space $ \H $ 
is a nondegenerate $ * $-homomorphism $ C^*_\max(G_q) \rightarrow \LH(\H) $. A basic example is the left regular representation 
of $ G_q $ on $ L^2(G_q)) $, which is constructed using the fundamental multiplicative unitary. 
The reduced group $ C^* $-algebra $ C^*_\red(G_q) $ of $ G_q $ is the image of $ C^*_\max(G_q) $ under the left regular representation. 

The definition of $ G_q $ given above makes perfect sense for $ q = 1 $, by starting with the classical algebra of functions $ C(K) $ and the 
group $ C^* $-algebra $ C^*(K) $ instead of their deformed versions. It is important to keep in mind, however, that the resulting function algebra $ C_0(G_1) $ 
is far from being isomorphic to the algebra of functions on the group $ G $. We also note that the group $ C^* $-algebra $ C^*(G_1) $ 
can be identified with the crossed product $ K \ltimes_\ad C(K) $ of the action of $ K $ on $ C(K) $ induced from the adjoint action of $ K $ on itself. 

In the proof of our main result we will need to consider two locally compact quantum groups closely related to $ G_1 $. Firstly, repeating the 
Drinfeld double construction with the maximal torus $ T $ instead of $ K $ yields the locally compact quantum group $ T \bowtie \hat{T} $, which is 
nothing but the direct product of the classical abelian groups $ T $ and $ \hat{T} $. 
Secondly, we will consider the relative Drinfeld double $ T \bowtie \hat{K} $, obtained by replacing $ C(K) $ with $ C(T) $ 
in the construction of $ C_0(G_1) $, while keeping the factor $ C^*(K) $. In this case the comultiplication involves 
the bicharacter $ (\pi \otimes \id)(W) $ instead of $ W $, obtained from the canonical restriction homomorphism $ \pi: C(K) \rightarrow C(T) $.

\section{Continuous fields of locally compact quantum groups} \label{seccflcqg} 

In this section we discuss continuous fields of locally compact quantum groups. Such continuous fields can be viewed as particular examples of locally 
compact quantum groupoids \cite{ELCQG}, for which many of the technical intricacies of the general theory of 
locally compact quantum groupoids disappear. Nonetheless, we will only sketch some key definitions and constructions. 

Let $ X $ be a locally compact Hausdorff space. Recall that a $ C_0(X) $-algebra is a $ C^* $-algebra $ A $ together with a 
nondegenerate $ * $-homomorphism $ C_0(X) \rightarrow ZM(A) $, where $ ZM(A) $ denotes the center of the multiplier algebra of $ A $. 
We will usually not distinguish between $ f \in C_0(X) $ and its image in $ ZM(A) $ in our notation, 
and write $ f a = a f $ for the product of $ f \in C_0(X) $ and $ a \in A $. 
If $ A $ and $ B $ are $ C_0(X) $-algebras then a $ * $-homomorphism $ \varphi: A \rightarrow B $ is called $ C_0(X) $-linear if $ \varphi(f a) = f \varphi(a) $ 
for all $ a \in A $ and $ f \in C_0(X) $. In this way $ C_0(X) $-algebras with $ C_0(X) $-linear $ * $-homomorphisms 
form naturally a category. The definition of $ C_0(X) $-linearity also applies to $ * $-homomorphisms $ \varphi: A \rightarrow M(B) $. 

Let $ I_x \subset C_0(X) $ be the kernel of evaluation at $ x \in X $, consisting of all $ f \in C_0(X) $ with $ f(x) = 0 $. 
The fiber of the $ C_0(X) $-algebra $ A $ at $ x $ is the $ C^* $-algebra $ A_x = A/[I_x A] $. The image of $ a \in A $ under the canonical 
quotient map $ A \rightarrow A_x $ will be denoted by $ a_x $, and we have 
$$
\|a \| = \sup_{x \in X} \|a_x\|
$$
for all $ a \in A $. Moreover, the norm function $ N_a: X \rightarrow \mathbb{R}, N_a(x) = \|a_x\| $ is upper semicontinuous and vanishes at infinity. 
A $ C_0(X) $-algebra $ A $ is called a continuous field of $ C^* $-algebras, or simply continuous, if $ N_a $ is continuous 
for all $ a \in A $. 
The trivial continuous field over $ X $ with fibre $ D $ is the $ C^* $-algebra $ A = C_0(X) \otimes D $ with the obvious $ C_0(X) $-algebra 
structure. Its fibers $ A_x $ identify canonically with $ D $ for all $ x \in X $. 

Let $ A $ be a $ C_0(X) $-algebra. A continuous field of representations of $ A $ on a Hilbert $ C_0(X) $-module $ \E $ is 
a $ C_0(X) $-linear $ * $-homomorphism $ \varphi: A \rightarrow \LH(\E) $, with respect to the canonical action $ C_0(X) \rightarrow M(\KH(\E)) = \LH(\E) $ 
obtained from the $ C_0(X) $-module structure on $ \E $. 
Such a continuous field induces representations $ \varphi_x: A_x \rightarrow \LH(\E_x) $ for all $ x \in X $, where $ \E_x = \E \otimes_{C_0(X)} \mathbb{C} $ 
is obtained by the taking interior tensor product with respect to evaluation at $ x \in X $. 
A continuous field of faithful representations is a continuous field of representations such that all maps $ \varphi_x $ are injective. 
If $ A $ admits a continuous field of faithful representations on some Hilbert $ C_0(X) $-module $ \E $ then $ A $ is continuous, and the converse holds 
if $ A $ is separable, see Th\'eor\`eme 3.3 in \cite{Blancharddef}. 

If $ A $ is a $ C_0(X) $-algebra then a $ C^* $-valued weight from $ A $ into $ C_b(X) = M(C_0(X)) $ is given by a hereditary cone $ P \subset A_+ $, 
closed under multiplication by elements from the positive part $ C_0(X)_+ $ of $ C_0(X) $,  and a positive linear map $ \omega $ from the linear 
span $ \mathfrak{M} $ of $ P $ into $ M(C_0(X)) $, such that $ \omega(f a) = f \omega(a) $ for all $ a \in \mathfrak{M} $ 
and $ f \in C_0(X)_+ $. We will only consider regular $ C^* $-valued weights in the sense of Kustermans, and refer to \cite{Kustermansregularweights} for 
more details.  

Let $ A, B $ be $ C_0(X) $-algebras. We define the $ C_0(X) $-tensor product $ A \otimes_{C_0(X)} B $ of $ A $ and $ B $ 
as in Definition 1.6 of \cite{Kasparov2}. If $ A $ and $ B $ are continuous fields of $ C^* $-algebras and $ X $ is compact then this 
agrees with the minimal balanced tensor product as defined by Blanchard in Definition 3.19 of \cite{Blancharddef}. 
The category of $ C_0(X) $-algebras becomes a monoidal category with this tensor product, with monoidal unit given by $ C_0(X) $ with 
its obvious $ C_0(X) $-algebra structure. 

Let us now define continuous fields of locally compact quantum groups, compare section 8.2 in \cite{ELCQG} and \cite{Enockcentralbasis}.  

\begin{definition} \label{defcflcqg}
Let $ X $ be a second countable locally compact space. A continuous field of locally compact quantum groups over $ X $ is a continuous field 
of $ C^* $-algebras $ C_0^\red({\bf G}) $ over $ X $ 
together with a nondegenerate $ * $-homomorphism $ \Delta: C_0^\red({\bf G}) \rightarrow M(C_0^\red({\bf G}) \otimes_{C_0(X)} C_0^\red({\bf G})) $ 
of $ C_0(X) $-algebras such that 
$$ 
(\Delta \otimes \id)\Delta = (\id \otimes \Delta)\Delta,  
$$ 
the induced maps $ \Delta_x: C_0^\red({\bf G})_x \rightarrow M(C_0^\red({\bf G})_x \otimes C_0^\red({\bf G})_x) $ are faithful for all $ x \in X $ and 
$$
[(C_0^\red({\bf G}) \otimes 1) \Delta(C_0^\red({\bf G}))] = C_0^\red({\bf G}) \otimes_{C_0(X)} C_0^\red({\bf G}) 
= [(1 \otimes C_0^\red({\bf G})) \Delta(C_0^\red({\bf G}))],
$$ 
together with $ C^* $-valued KMS-weights $ \phi $ and $ \psi $ from $ C_0^\red({\bf G}) $ into $ M(C_0(X)) $ which induce faithful left and right invariant 
KMS-weights on all fibres, respectively.  
\end{definition} 

We point out here that the fibers $ C_0^\red({\bf G})_x $ together with 
the comultiplications $ \Delta_x: C_0^\red({\bf G})_x \rightarrow M(C_0^\red({\bf G})_x \otimes C_0^\red({\bf G})_x) $ induced from $ \Delta $ are 
naturally Hopf $ C^* $-algebras. 
For the definition of $ C^* $-values KMS-weights we refer to \cite{ELCQG} and references therein. The requirements on $ \phi $ and $ \psi $ 
in Definition \ref{defcflcqg} ensure that the induced weights $ \phi_x, \psi_x $ 
on $ C_0^\red({\bf G})_x $ are left and right invariant faithful KMS-weights in the sense of \cite{KVLCQG} for all $ x \in X $. 
In particular, the fibers $ C_0^\red({\bf G})_x $ of a continuous field in the sense of Definition \ref{defcflcqg} 
determine locally compact quantum groups, which we will denote by $ {\bf G}_x $ in the sequel. 
We will assume throughout that the underlying $ C^* $-algebras of our continuous fields of locally compact quantum groups are separable. 

From the $ C^* $-valued weight $ \phi $ one obtains a Hilbert $ C_0(X) $-module $ L^2({\bf G}) $ and a continuous field of 
multiplicative unitaries $ {\bf W} $ on $ L^2({\bf G}) \otimes_{C_0(X)} L^2({\bf G}) $, compare Definition 4.6 in \cite{Blancharddef}, 
such that $ {\bf W}_x $ is the fundamental multiplicative unitary of $ {\bf G}_x $ for all $ x \in X $. 
Moreover, the comultiplication $ \Delta $ can be written as $ \Delta(f) = {\bf W}^*(1 \otimes f) {\bf W} $ 
as in the case of locally compact quantum groups. 
Let us say that a continuous field of locally compact quantum groups is regular if the corresponding continuous field of multiplicative unitaries $ {\bf W} $ 
is regular, see section 4.2 in \cite{Blancharddef}. We will only be interested in regular continuous fields in the sequel. 

The dual of a continuous field of locally compact quantum groups $ {\bf G} $ is obtained from the right leg of $ {\bf W} $ by setting 
$$
C^*_\red({\bf G}) = [(\omega \otimes \id)({\bf W}) \mid \omega \in \LH(L^2({\bf G}))_* ],
$$
where $ \LH(L^2({\bf G}))_* $ is the continuous field of preduals as in D\'efinition 4.2 of \cite{Blancharddef}. 
Together with the comultiplication $ \hat{\Delta}(x) = \hat{{\bf W}}^*(1 \otimes x) \hat{{\bf W}} $ for $ \hat{{\bf W}} = \Sigma {\bf W}^* \Sigma $ 
this can be made again a continuous field of locally compact quantum groups. 
We will write $ \hat{\bf G} $ for the resulting continuous field, so that we have $ C_0^\red(\hat{\bf G}) = C^*_\red({\bf G}) $ with 
comultiplication $ \hat{\Delta} $, and we write $ \check{\bf G} $ for the continuous field obtained by swapping the comultiplication in $ C_0^\red(\hat{\bf G}) $. 

Let us now discuss the key example of a continuous field of locally compact quantum groups that we will be interested in. 
We fix $ q = e^h > 0 $ and consider a simply connected semisimple Lie group $ G $ with maximal compact subgroup $ K $. 
With the notations as in section \ref{secprelim}, there are natural continuous fields of $ C^* $-algebras $ C({\bf K}) $ and $ C^*({\bf K}) $ over $ [0,1] $ 
with fibers $ C({\bf K})_\sigma = C(K_{q^\sigma}) $ and $ C^*({\bf K})_\sigma = C^*(K_{q^\sigma}) $,
which assemble to a continuous field of $ C^* $-algebras 
$$ 
C_0({\bf G}) = C({\bf K}) \otimes_{C[0,1]} C^*({\bf K}) 
$$ 
with fibers $ C_0({\bf G})_\sigma = C_0(G_{q^\sigma}) $, see \cite{Blancharddef}, \cite{MVcqgdeformation}. 
The fiberwise coproducts determine a coproduct $ \Delta_{\bf G}: C_0({\bf G}) \rightarrow M(C_0({\bf G}) \otimes_{C[0,1]} C_0({\bf G})) $ 
which turns $ C_0({\bf G}) $ into a regular continuous field of locally compact quantum groups over $ [0,1] $.

\section{Equivariant Kasparov theory} \label{seckkg} 

In this section we develop the basics of equivariant $ KK $-theory for continuous fields of locally compact quantum groups. 
Throughout let $ X $ be a fixed second countable locally compact space.  

Let $ {\bf G} $ be a continuous field of locally compact quantum groups over $ X $, and let $ A $ be a $ C_0(X) $-algebra. 
An action of $ {\bf G} $ on $ A $ is an injective nondegenerate $ C_0(X) $-linear $ * $-homomorphism $ \alpha: A \rightarrow M(C_0^\red({\bf G}) \otimes_{C_0(X)} A) $ 
such that $ (\Delta \otimes \id) \alpha = (\id \otimes \alpha) \alpha $ and $ [(C_0^\red({\bf G}) \otimes 1) \alpha(A)] = C_0^\red({\bf G}) \otimes_{C_0(X)} A $. 
We will also say that $ A $ is a $ {\bf G} $-$ C^* $-algebra in this case. 
A $ {\bf G} $-equivariant $ * $-homomorphism between $ {\bf G} $-$ C^* $-algebras $ A, B $ with coactions $ \alpha, \beta $, respectively, 
is a $ C_0(X) $-linear $ * $-homomorphism $ \varphi: A \rightarrow B $ such that $ \beta \varphi = (\id \otimes \varphi) \alpha $. 
This notion of equivariance makes sense also for $ * $-homomorphisms $ \varphi: A \rightarrow M(B) $. 
The trivial action of $ {\bf G} $ is given by $ C_0(X) $, viewed as a $ C_0(X) $-algebra over itself, with the canonical 
map $ C_0(X) \rightarrow M(C_0^\red({\bf G}) \otimes_{C_0(X)} C_0(X)) \cong M(C_0^\red({\bf G})) $ as coaction. 

Let $ {\bf G} $ be a continuous field of locally compact quantum groups and let $ \beta: B \rightarrow M(C_0^\red({\bf G}) \otimes_{C_0(X)} B) $ be an action 
of $ {\bf G} $ on the $ C^* $-algebra $ B $. An action of $ {\bf G} $ on a Hilbert $ B $-module $ \E $ is a nondegenerate $ C_0(X) $-linear 
morphism $ \lambda: \E \rightarrow M(C_0^\red({\bf G}) \otimes_{C_0(X)} \E) $ such that $ (\Delta \otimes \id) \lambda = (\id \otimes \lambda)\lambda $ and 
$$ 
[(C_0^\red({\bf G}) \otimes 1)\lambda(\E)] = C_0^\red({\bf G}) \otimes_{C_0(X)} \E = [\lambda(\E)(C_0^\red({\bf G}) \otimes 1)];
$$
compare \cite{NVpoincare} for the notion of multiplier modules of Hilbert modules and their morphisms. We will also say that $ \E $ is 
a $ {\bf G} $-Hilbert $ B $-module in this case. A $ {\bf G} $-Hilbert space is a $ {\bf G} $-Hilbert module over 
the trivial $ {\bf G} $-$ C^* $-algebra $ B = C_0(X) $. 

If $ \E $ is a $ {\bf G} $-Hilbert $ B $-module with coaction $ \lambda: \E \rightarrow M(C_0^\red({\bf G}) \otimes_{C_0(X)} \E) $ 
then we obtain a unitary $ V_\lambda: \E \otimes_B (C_0^\red({\bf G}) \otimes_{C_0(X)} B) \rightarrow C_0^\red({\bf G}) \otimes_{C_0(X)} \E $
satisfying $ V_\lambda(\xi \otimes x) = \lambda(\xi) x $ for $ \xi \in \E, x \in C_0^\red({\bf G}) \otimes_{C_0(X)} B $. 
Moreover $ \KH(\E) $ becomes a $ {\bf G} $-$ C^* $-algebra with the adjoint 
action $ \ad_\lambda: \KH(\E) \rightarrow M(C_0^\red({\bf G}) \otimes_{C_0(X)} \KH(\E)) \cong \LH(C_0^\red({\bf G}) \otimes_{C_0(X)} \E) $ given by 
$$
\ad_\lambda(T) = V_\lambda(T \otimes \id) V_\lambda^*.  
$$

Let us next discuss crossed products and duality. If $ A $ is a $ {\bf G} $-$ C^* $-algebra then the reduced crossed product $ {\bf G} \ltimes_\red A $ is 
defined by 
$$
{\bf G} \ltimes_\red A = [(C^*_\red({\bf G}) \otimes 1) \alpha(A)] \subset \LH(L^2({\bf G}) \otimes_{C_0(X)} A). 
$$
Note that the reduced crossed product is naturally a $ C_0(X) $-algebra. 
The dual action $ \check{\alpha}: {\bf G} \ltimes_\red A \rightarrow M(C_0^\red(\check{\bf G}) \otimes_{C_0(X)} {\bf G} \ltimes_\red A) $ 
on the reduced crossed product is given by $ \check{\alpha}(x) = \check{\bf V}_{12}(1 \otimes x) \check{\bf V}_{12} $, compare 
Theorem 4.2 in \cite{Timmermannhopfcstar}. Here $ \check{\bf V} = \Sigma \hat{\bf V} \Sigma $, where $ \hat{\bf V} $ 
is the continuous field of multiplicative unitaries constructed from $ {\bf W} $ in the same way as for locally compact 
quantum groups, see section \ref{secprelim}. 

Let $ \KH_{\bf G} = \KH(L^2({\bf G})) $ and consider $ \KH_{\bf G} \otimes_{C_0(X)} A \cong \KH(L^2({\bf G}) \otimes_{C_0(X)} A) $ 
as a $ {\bf G} $-$ C^* $-algebra with the conjugation action of the action $ \lambda(x \otimes a) = {\bf X}_{12}^* \Sigma_{12} (\id \otimes \alpha)(x \otimes a) $ 
on $ L^2({\bf G}) \otimes_{C_0(X)} A $, where $ {\bf X} = \Sigma {\bf V} \Sigma $, compare \cite{NVpoincare}. 
As a special case of the results in \cite{Timmermannhopfcstar} we obtain the following version of Takesaki-Takai duality.  

\begin{theorem} \label{TTduality}
Let $ {\bf G} $ be a regular continuous field of locally compact quantum groups and let $ A $ be a $ {\bf G} $-$ C^* $-algebra. 
Then $ \check{\bf G} \ltimes_\red {\bf G} \ltimes_\red A $ is naturally $ {\bf G} $-equivariantly isomorphic to $ \KH_{\bf G} \otimes_{C_0(X)} A $. 
\end{theorem} 

We will now sketch the construction of equivariant $ KK $-theory for continuous fields in the sense of Definition \ref{defcflcqg}, in analogy to the case of 
locally compact quantum groups \cite{BSKK}, \cite{NVpoincare}. 
Let $ {\bf G}  $ be a continuous field of locally compact quantum groups over $ X $. 
Moreover let $ A $ and $ B $ be separable $ {\bf G} $-$ C^* $-algebras. 
A $ {\bf G} $-equivariant Kasparov $ A $-$ B $-module is a countably generated graded $ {\bf G} $-equivariant Hilbert $ B $-module $ \E $ together 
with a $ {\bf G} $-equivariant $ * $-homomorphism $ \phi: A \rightarrow \LH(\E) $ and an odd operator $ F \in \LH(\E) $ such that
$$
[F, \phi(a)], \qquad (F^2 - 1) \phi(a), \qquad (F - F^*)\phi(a)
$$
are contained in $ \KH(\E) $ for all $ a \in A $, and $ F $ is almost invariant in the sense that
$$
(\id \otimes \phi)(x)(1 \otimes F - \ad_\lambda(F)) \in C_0^\red({\bf G}) \otimes_{C_0(X)} \KH(\E)
$$
for all $ x \in C_0^\red({\bf G}) \otimes_{C_0(X)} A $. 
Here $ C_0^\red({\bf G}) \otimes_{C_0(X)} \KH(\E) = \KH(C_0^\red({\bf G}) \otimes_{C_0(X)} \E) $ 
is viewed as a subset of $ \LH(C_0^\red({\bf G}) \otimes_{C_0(X)} \E) $ and $ \ad_\lambda $ is the adjoint action associated to the given 
action $ \lambda: \E \rightarrow M(C_0^\red({\bf G}) \otimes_{C_0(X)} \E) $ on $ \E $.
Two $ {\bf G} $-equivariant Kasparov $ A $-$ B $-modules $ (\E_0, \phi_0, F_0) $ and $ (\E_1, \phi_1, F_1) $ are called unitarily equivalent if there is
a $ {\bf G} $-invariant unitary $ U \in \LH(\E_0, \E_1) $ of degree zero such that $ U \phi_0(a) = \phi_1(a) U $ for all $ a \in A $
and $ F_1 U = U F_0 $. We write $ (\E_0, \phi_0, F_0) \cong (\E_1, \phi_1, F_1) $ in this case.
Let $ E_{\bf G}(A,B) $ be the set of unitary equivalence classes of $ {\bf G} $-equivariant Kasparov $ A $-$ B $-modules.
This set is functorial for $ {\bf G} $-equivariant $ * $-homomorphisms
in both variables. 
A homotopy between $ {\bf G} $-equivariant Kasparov $ A $-$ B $-modules $ (\E_0, \phi_0, F_0) $ and $ (\E_1, \phi_1, F_1) $ is a
$ {\bf G} $-equivariant Kasparov $ A $-$ B[0,1] $-module
$ (\E, \phi, F) $ such that $ (\ev_t)_*(\E, \phi, F) \cong (\E_t, \phi_t, F_t) $ for $ t = 0,1 $. Here
$ B[0,1] = B \otimes C[0,1] $ is equipped with the action induced from $ B $ 
and $ \ev_t: B[0,1] \rightarrow B $ is evaluation at $ t $.

\begin{definition} \label{defkkg}
Let $ {\bf G} $ be a continuous field of locally compact quantum groups and let $ A $ and $ B $ be $ {\bf G} $-$ C^* $-algebras. 
The $ {\bf G} $-equivariant Kasparov group $ KK^{\bf G}(A,B) $ is the set of homotopy classes of $ {\bf G} $-equivariant Kasparov $ A $-$ B $-modules.
\end{definition}

We note that $ KK^{\bf G}(A,B) $ becomes an abelian group with addition given by the direct sum of Kasparov modules.
Many properties of ordinary $ KK $-theory and equivariant $ KK $-theory for locally compact quantum groups carry over to 
the $ {\bf G} $-equivariant situation. This includes in particular the construction of the
Kasparov composition product and Bott periodicity \cite{BSKK}, but we will not spell out the details. 
As usual we write $ KK^{\bf G}_0(A,B) = KK^{\bf G}(A,B) $ and let $ KK^{\bf G}_1(A,B) $ be the odd $ KK $-group obtained by suspension in either variable.
We note that in the case that $ {\bf G} $ is the trivial field of trivial (quantum) groups over $ X $, Definition \ref{defkkg} 
reduces to $ \R KK(X; A,B) $ in the sense of Kasparov \cite{Kasparov2}. 

Using Theorem \ref{TTduality} one obtains the following version of the Baaj-Skandalis duality theorem \cite{BSKK}. 

\begin{theorem}\label{BSduality}
Let $ {\bf G} $ be a regular continuous field of locally compact quantum groups. 
For all $ {\bf G} $-$ C^* $-algebras $ A $ and $ B $ there is a canonical isomorphism
$$
J_{\bf G}: KK^{\bf G}(A,B) \rightarrow KK^{\check{\bf G}}({\bf G} \ltimes_\red A, {\bf G} \ltimes_\red B), 
$$
which is multiplicative with respect to the Kasparov product.
\end{theorem}

We note that under this isomorphism the class of a $ {\bf G} $-equivariant Kasparov $ A $-$ B $-module $ (\E, \phi, F) $ is mapped
to the class of a $ \check{\bf G} $-equivariant Kasparov module $ (J_{\bf G}(\E), J_{\bf G}(\phi), J_{\bf G}(F)) $ with an operator $ J_{\bf G}(F) $ 
which is exactly invariant under the coaction of $ C_0^\red(\hat{\bf G}) $. 

Let $ {\bf G} $ be a regular continuous field of locally compact quantum groups and let $ \E $ and $ \F $ be $ {\bf G} $-Hilbert $ B $-modules 
which are isomorphic as Hilbert $ B $-modules. Then we have a $ {\bf G} $-equivariant isomorphism
$$
L^2({\bf G}) \otimes_{C_0(X)} \E \cong L^2({\bf G}) \otimes_{C_0(X)} \F
$$
of $ {\bf G} $-Hilbert $ B $-modules where $ L^2({\bf G}) $ is viewed as a $ {\bf G} $-Hilbert space using the left regular 
representation. Using the Kasparov stabilisation theorem we deduce the equivariant stabilisation theorem, namely that there is a $ {\bf G} $-equivariant 
Hilbert $ B $-module isomorphism
$$
(L^2({\bf G}) \otimes_{C_0(X)} \E) \oplus (L^2({\bf G}) \otimes_{C_0(X)} (\HH \otimes B)) \cong L^2({\bf G}) \otimes_{C_0(X)} (\HH \otimes B)
$$
for every countably generated $ {\bf G} $-Hilbert $ B $-module $ \E $, where $ \HH $ is some fixed infinite dimensional separable Hilbert space. 

It follows from Baaj-Skandalis duality that $ KK^{\bf G}(A,B) $ can be represented by homotopy classes of $ {\bf G} $-equivariant Kasparov
$ (\KH_{\bf G} \otimes_{C_0(X)} A) $-$ (\KH_{\bf G} \otimes_{C_0(X)} B) $-modules $ (\E, \phi, F) $ with $ {\bf G} $-invariant operator $ F $. 
Using the $ {\bf G} $-equivariant Morita equivalence between $ \KH_{\bf G} \otimes_{C_0(X)} B $ and $ B $ we see that $ KK^{\bf G}(A,B) $ can be 
represented by homotopy classes of equivariant Kasparov $ (\KH_{\bf G} \otimes_{C_0(X)} A) $-$ B $ modules of the form $ (L^2({\bf G}) \otimes_{C_0(X)} \E, \phi, F) $ 
with invariant $ F $. Using equivariant stabilisation we can furthermore assume 
that $ (L^2({\bf G}) \otimes_{C_0(X)} \E)_\pm = L^2({\bf G}) \otimes_{C_0(X)} (\HH \otimes B) $ is the standard $ {\bf G} $-Hilbert $ B $-module. 

Our next aim is to establish the Cuntz picture of $ KK^{\bf G} $ for regular continuous fields of locally compact quantum groups.
This can be done in a similar way as in \cite{Meyerkkg}, \cite{NVpoincare}. 
Let $ {\bf G} $ be a continuous field of locally compact quantum groups and let $ A_1 $ and $ A_2 $ be $ {\bf G} $-$ C^* $-algebras. 
In analogy to the construction of $ A_1 \otimes_{C_0(X)} A_2 $ in \cite{Kasparov2} 
we define the $ C_0(X) $-free product $ A_1 *_{C_0(X)} A_2 $ as completion of the image of the algebraic free product of $ A_1 $ and $ A_2 $ 
in all possible common $ C_0(X) $-representations of $ A_1, A_2 $ on Hilbert modules over the bidual $ C_0(X)^{**} $. 
It is straightforward to check that $ A_1 *_{C_0(X)} A_2 $ is a $ C_0(X) $-algebra such that the 
canonical $ * $-homomorphisms $ \iota_j: A_j \rightarrow A_1 *_{C_0(X)} A_2 $ for $ j = 1,2 $ are $ C_0(X) $-linear. 
Using the universal property of $ A_1 *_{C_0(X)} A_2 $ 
we obtain a $ C_0(X) $-linear $ * $-homomorphism $ \alpha: A_1 *_{C_0(X)} A_2 \rightarrow M(C_0^\red({\bf G}) \otimes_{C_0(X)} (A_1 *_{C_0(X)} A_2)) $. This map
satisfies all properties of an action except that it is not
clear whether $ \alpha $ is always injective. If necessary, this technicality can be arranged by passing to the quotient
of $ A_1 *_{C_0(X)} A_2 $ by the kernel of $ \alpha $. By abuse of notation we will write again $ A_1 *_{C_0(X)} A_2 $ for the 
resulting $ {\bf G} $-$ C^* $-algebra, although this quotient may depend on $ {\bf G} $. 
The resulting $ {\bf G} $-$ C^* $-algebra is universal for pairs of $ {\bf G} $-equivariant $ * $-homomorphisms $ f_1: A_1 \rightarrow C $ and
$ f_2 : A_2 \rightarrow C $ into $ {\bf G} $-$ C^* $-algebras $ C $. That is, for any such pair of $ * $-homomorphisms there exists a unique
$ {\bf G} $-equivariant $ * $-homomorphism $ f: A_1 *_{C_0(X)} A_2 \rightarrow C $ such that $ f \iota_j = f_j $ for $ j = 1,2 $.

Let $ A $ be an $ {\bf G} $-$ C^* $-algebra and consider $ QA = A *_{C_0(X)} A $. There is a canonical extension
\begin{equation*}
\xymatrix{
0 \ar@{->}[r] & qA \ar@{->}[r] & QA \ar@{->}[r]^{\pi} & A \ar@{->}[r] & 0
}
\end{equation*}
of $ {\bf G} $-$ C^* $-algebras with $ {\bf G} $-equivariant splitting; here $ \pi $ is the homomorphism associated to the pair $ f_1 = \id_A = f_2 $ and
$ qA $ its kernel. 

Writing $ [A,B]_{\bf G} $ for the set of equivariant homotopy classes of $ {\bf G} $-equivariant $ * $-homomorphisms 
between $ {\bf G} $-$ C^* $-algebras $ A $ and $ B $ and $ \KH $ for the compact operators on some separable Hilbert space, 
we arrive at the following description of the equivariant $ KK $-groups.

\begin{theorem} \label{cuntzpic}
Let $ {\bf G} $ be a regular continuous field of locally compact quantum groups. Then there is a natural isomorphism
$$
KK^{\bf G}(A,B) \cong [q(\KH_{\bf G} \otimes_{C_0(X)} A), \KH_{\bf G} \otimes_{C_0(X)} (\KH \otimes B)]_{\bf G}
$$
for all separable $ {\bf G} $-$ C^* $-algebras $ A $ and $ B $. If we write $ \KH(\bf G) = \KH_{\bf G} \otimes \KH $ we also have a natural isomorphism
$$
KK^{\bf G}(A,B) \cong [\KH({\bf G}) \otimes_{C_0(X)} q(\KH({\bf G}) \otimes_{C_0(X)} A), 
\KH({\bf G}) \otimes_{C_0(X)} q(\KH({\bf G}) \otimes_{C_0(X)} B)]_{\bf G}
$$
under which the Kasparov product corresponds to the composition of homomorphisms.
\end{theorem}

Consider the category $ {\bf G}\hyphenAlg $ of all separable $ {\bf G} $-$ C^* $-algebras and $ {\bf G} $-equivariant $ * $-homomorphisms 
for a regular continuous field of locally compact quantum groups $ {\bf G} $. A functor $ F $ 
from $ {\bf G}\hyphenAlg $ to an additive category $ \C $ is called a homotopy functor if $ F(f_0) = F(f_1) $ whenever $ f_0 $ and $ f_1 $ 
are $ {\bf G} $-equivariantly homotopic. It is called stable if for all separable $ {\bf G} $-Hilbert spaces $ \H_1, \H_2 $
the maps $ F(\KH(\H_j) \otimes_{C_0(X)} A) \rightarrow F(\KH(\H_1 \oplus \H_2) \otimes_{C_0(X)} A) $
induced by the canonical inclusions $ \H_j \rightarrow \H_1 \oplus \H_2 $ for $ j = 1,2 $
are isomorphisms. As in the group case, a homotopy functor
$ F $ is stable iff there exists a natural isomorphism
$ F(A) \cong F(\KH({\bf G}) \otimes_{C_0(X)} A) $ for all $ A $.
Finally, $ F $ is called split exact if for every extension
\begin{equation*}
\xymatrix{
0 \ar@{->}[r] & K \ar@{->}[r] & E \ar@{->}[r] & Q \ar@{->}[r] & 0
}
\end{equation*}
of $ {\bf G} $-$ C^* $-algebras that splits by a $ {\bf G} $-equivariant $ * $-homomorphism $ \sigma: Q \rightarrow E $
the induced sequence $ 0 \rightarrow F(K) \rightarrow F(E) \rightarrow F(Q) \rightarrow 0 $ in $ \C $ is split exact. 

Equivariant $ KK $-theory can be viewed as an additive category $ KK^{\bf G} $ with separable
$ {\bf G} $-$ C^* $-algebras as objects and $ KK^{\bf G}(A,B) $ as the set of morphisms between objects $ A $ and $ B $.
Composition of morphisms is given by the Kasparov product.
There is a canonical functor $ \iota: {\bf G}\hyphenAlg \rightarrow KK^{\bf G} $ which is the identity on
objects and sends equivariant $ * $-homomorphisms to the corresponding $ KK $-elements. This functor is a split exact
stable homotopy functor. 
As a consequence of theorem \ref{cuntzpic} we obtain the following universal property of $ KK^{\bf G} $, compare again \cite{Meyerkkg}.

\begin{theorem} \label{kkuniversal}
Let $ {\bf G} $ be a regular continuous field of locally compact quantum groups. The functor $ \iota: {\bf G} \hyphenAlg \rightarrow KK^{\bf G} $ 
is the universal split exact stable homotopy functor on the category $ {\bf G} \hyphenAlg $.
More precisely, if $ F: {\bf G} \hyphenAlg \rightarrow \C  $ is any split exact stable homotopy functor
with values in an additive category $ \C $ then there exists a unique functor $ f: KK^{\bf G} \rightarrow \C $ such that
$ F = f \iota $.
\end{theorem}

Let us also explain how $ KK^{\bf G} $ becomes a triangulated category, compare \cite{MNtriangulated}. 
Let $ \Sigma A $ denote the suspension $ C_0(\mathbb{R}) \otimes A $ of a $ {\bf G} $-$ C^* $-algebra $ A $.
Here $ C_0(\mathbb{R}) $ is equipped with the trivial coaction. The corresponding functor $ \Sigma: KK^{\bf G} \rightarrow KK^{\bf G} $ determines the
translation automorphism.
If $ \varphi: A \rightarrow B $ is a $ {\bf G} $-equivariant $ * $-homomorphism then the mapping cone
$$
C_\varphi = \{(a,b) \in A \times C_0((0,1], B) | b(1) = \varphi(a) \}
$$
is a $ {\bf G} $-$ C^* $-algebra in a natural way, and there is a canonical diagram
$$
\xymatrix{
\Sigma B  \;\; \ar@{->}[r] & C_\varphi \ar@{->}[r] & A \ar@{->}[r]^\varphi & B
}
$$
of $ {\bf G} $-equivariant $ * $-homomorphisms. Diagrams of this form are called mapping cone triangles.
By definition, an exact triangle is a diagram $ \Sigma Q \rightarrow K \rightarrow E \rightarrow Q $ in $ KK^{\bf G} $
which is isomorphic to a mapping cone triangle. 

The proof of the following result is carried out in the same way as for regular locally compact quantum 
groups \cite{MNtriangulated}, \cite{NVpoincare}.

\begin{prop}\label{kktriang}
Let $ {\bf G} $ be a regular continuous field of locally compact quantum groups. Then the category $ KK^{\bf G} $ together with the translation functor
and the exact triangles described above is triangulated.
\end{prop}

We will need a version of the Green-Julg theorem for continuous fields of compact quantum groups. Here a continuous field of locally 
compact quantum groups $ {\bf G} $ is called compact if all fibers of $ C_0^\red({\bf G}) $ are compact quantum groups. 
If $ A $ is a $ C_0(X) $-algebra we write $ \res^E_{\bf G}(A) $ for the $ {\bf G} $-$ C^* $-algebra $ A $ with the trivial action, 
given by sending $ a \in A $ to $ 1 \otimes a \in M(C_0^\red({\bf G}) \otimes_{C_0(X)} A) $. 

\begin{theorem} \label{GJ}
Let $ {\bf G} $ be a continuous field of compact quantum groups. Then there is a natural isomorphism
$$
KK^{\bf G}(\res^E_{\bf G}(A), B) \cong \R KK(X; A, {\bf G} \ltimes_\red B)
$$
for all $ C_0(X) $-algebras $ A $ and all $ {\bf G} $-$ C^* $-algebras $ B $. 
\end{theorem}

The proof is carried out in the same way as for compact quantum groups, using parametrised versions of the construction of the unit and counit of the 
adjunction. 

If $ {\bf G} $ is a continuous field of locally compact quantum groups over $ X $, then a field of open quantum subgroups is a continuous field 
of locally compact quantum groups $ {\bf H} $ over $ X $ together with a $ C_0(X) $-linear $ * $-homomorphism $ \pi: C^\red_0({\bf G}) \rightarrow C^\red_0({\bf H}) $ 
compatible with the comultiplications, exhibiting all fibers $ {\bf H}_x $ as open quantum subgroups of $ {\bf G}_x $. 
The map $ \pi $ is induced from a central projection $ 1_{\bf H} \in M(C^\red_0({\bf G})) $ and 
the $ C^* $-valued invariant weights of $ C^\red_0({\bf H}) $ are obtained by restriction from those of $ C^\red_0({\bf G}) $. 

Assume that $ B $ is a $ {\bf G} $-$ C^* $-algebra with coaction $ \beta: B \rightarrow M(C_0^\red({\bf G}) \otimes_{C_0(X)} B) $. Then 
the restriction $ \res^{\bf G}_{\bf H}(B) $ is the $ {\bf H} $-$ C^* $-algebra defined 
by $ (\pi \otimes \id) \beta: B \rightarrow M(C_0^\red({\bf H}) \otimes_{C_0(X)} B) $. 
Using Theorem \ref{kkuniversal} one checks that restriction lifts to a functor $ \res^{\bf G}_{\bf H}: KK^{\bf G} \rightarrow KK^{\bf H} $. 

Conversely, following our discussion in section \ref{secprelim}, let us sketch how to define an induced $ {\bf G} $-$ C^* $-algebra $ \ind_{\bf H}^{\bf G}(A) $ 
for any $ {\bf H} $-$ C^* $-algebra $ A $ with coaction $ \alpha: A \rightarrow M(C_0^\red({\bf H}) \otimes_{C_0(X)} A) $, provided the continuous field 
of quantum groups $ {\bf G} $ is regular and $ {\bf H} \subset {\bf G} $ is open. 
In fact, with the obvious modifications, the construction of the imprimitivity bimodule $ \J = [(\I \otimes 1)\alpha(A)] $ 
discussed in section \ref{secprelim} makes sense in this setting, and allows us to define 
a $ {\bf G} $-$ C^* $-algebra $ \ind_{\bf H}^{\bf G}(A) $, 
such that $ \J $ defines a covariant Morita equivalence $ {\bf G} \ltimes_\red \ind_{\bf H}^{\bf G}(A) \sim_M {\bf H} \ltimes_\red A $ 
of $ \check{\bf G} $-$ C^* $-algebras. 
Moreover, there is a projection $ p = 1_{\bf H} \otimes \id \in M({\bf G} \ltimes_\red \ind_{\bf H}^{\bf G}(A)) $ such 
that $ p $ exhibits $ {\bf H} \ltimes_\red A $ as a full corner of $ {\bf G} \ltimes_\red \ind_{\bf H}^{\bf G}(A) $. 
The construction of induced algebras lifts to a functor $ \ind_{\bf H}^{\bf G}: KK^{\bf H} \rightarrow KK^{\bf G} $ on the level of $ KK $-theory. 

The following result is a generalisation of Proposition 6.2 in \cite{VVfreeu}. 

\begin{theorem} \label{opensubgroupinduction}
Let $ {\bf G} $ be a regular continuous field of locally compact quantum groups and let $ {\bf H} \subset {\bf G} $ be a continuous 
family of open quantum subgroups of $ {\bf G} $. Then there exists an isomorphism 
$$
KK^{\bf G}_*(\ind_{\bf H}^{\bf G}(A), B) \cong KK^{\bf H}_*(A, \res^{\bf G}_{\bf H}(B)),  
$$
natural in $ A \in KK^{\bf H}, B \in KK^{\bf G} $.  
\end{theorem} 

\proof For $ A \in KK^{\bf H} $ there is a natural $ \check{\bf G} $-equivariant 
$ * $-homomorphism $ {\bf H} \ltimes_\red A \rightarrow {\bf H} \ltimes_\red \res^{\bf G}_{\bf H} \ind_{\bf H}^{\bf G}(A) $ 
coming from the embedding 
$$ 
[(C^*_\red({\bf H}) \otimes 1) \alpha(A)] \cong [(i \otimes \id)(C^*_\red({\bf H}) \otimes 1) \alpha(A)(i \otimes \id)^*] \subset [\J \J^*], 
$$
where $ i: L^2({\bf H}) \rightarrow L^2({\bf G}) $ is the canonical inclusion. 
By Baaj-Skandalis duality, this homomorphism is of the form $ {\bf H} \ltimes_\red \eta_A $ 
for $ \eta_A: A \rightarrow \res^{\bf G}_{\bf H} \ind_{\bf H}^{\bf G}(A) $ in $ KK^{\bf H} $. 

Similarly, for $ B \in KK^{\bf G} $ we obtain a 
morphism $ {\bf G} \ltimes_\red \ind_{\bf H}^{\bf G} \res^{\bf G}_{\bf H}(B) \rightarrow {\bf G} \ltimes_\red B $ 
in $ KK^{\check{\bf G}} $ 
by composing the Morita equivalence $ {\bf G} \ltimes_\red \ind_{\bf H}^{\bf G} \res^{\bf G}_{\bf H}(B) \sim_M {\bf H} \ltimes_\red \res^{\bf G}_{\bf H}(B) $ 
with the canonical inclusion $ {\bf H} \ltimes_\red \res^{\bf G}_{\bf H}(B) \rightarrow {\bf G} \ltimes_\red B $. 
This is of the form $ {\bf G} \ltimes_\red \kappa_B $ for $ \kappa_B: \ind_{\bf H}^{\bf G} \res^{\bf G}_{\bf H}(B) \rightarrow B $ in $ KK^{\bf G} $. 

Let us check that these morphisms satisfy the required identities for an adjunction. The composition 
$$
\xymatrix{
{\bf G} \ltimes_\red \ind_{\bf H}^{\bf G}(A) \ar@{->}[r]^{\!\!\!\!\!\!\!\!\!\!\!\!{\bf G} \ltimes_\red \ind(\eta_A)} \;\;\; & \;\;\;\; 
{\bf G} \ltimes_\red \ind_{\bf H}^{\bf G} \res^{\bf G}_{\bf H} \ind_{\bf H}^{\bf G}(A) 
\ar@{->}[r]^{\qquad \;\; {\bf G} \ltimes_\red \kappa_{\ind(A)}} \;\;\; & \;{\bf G} \ltimes_\red \ind_{\bf H}^{\bf G}(A)
}
$$
in $ KK^{\check{\bf G}} $ can be identified with 
$$
\xymatrix{
{\bf H} \ltimes_\red A \ar@{->}[r]^{\!\!\!\!\!\!\!\!\!\!\!\!\!\!\!\!\!\!\!\! {\bf H} \ltimes_\red \eta_A} 
& \;\; {\bf H} \ltimes_\red \res^{\bf G}_{\bf H} \ind_{\bf H}^{\bf  G}(A) 
\ar@{->}[r]^{} & \;{\bf G} \ltimes_\red \ind_{\bf H}^{\bf G}(A) \sim_M {\bf H} \ltimes_\red A.
}
$$
The map $ {\bf H} \ltimes_\red A \rightarrow {\bf G} \ltimes_\red \ind_{\bf H}^{\bf G}(A) $ in here is the canonical inclusion. 
It follows that the composition $ \kappa_{\ind A} \circ \ind(\eta_A) $ is the identity. 

Similarly, the composition 
$$
\xymatrix{
{\bf H} \ltimes_\red \res^{\bf G}_{\bf H}(B) \ar@{->}[r]^{\!\!\!\!\!\!\!\!\!\!\!\!\!\!\!\!\!\! {\bf H} \ltimes_\red \eta_{\res(B)}} \;
& \;\;\;\; {\bf H} \ltimes_\red \res^{\bf G}_{\bf H} \ind_{\bf H}^{\bf G} \res^{\bf G}_{\bf H}(B) 
\ar@{->}[r]^{\qquad \;\;\; {\bf H} \ltimes_\red \res(\kappa_B)} \;\; & \;\; {\bf H} \ltimes_\red \res^{\bf G}_{\bf H}(B)
}
$$
in $ KK^{\check{\bf H}} $ corresponds to 
$$
\xymatrix{
{\bf H} \ltimes_\red \res^{\bf G}_{\bf H}(B) \ar@{->}[r]^{\!\!\!\!\!\!\!\!\!\!\! {\bf H} \ltimes_\red \eta_{\res(B)} \;\;\;\;} \;\;
& \;\;\;\; {\bf H} \ltimes_\red \res^{\bf G}_{\bf H} \ind_{\bf H}^{\bf G} \res^{\bf G}_{\bf H}(B) 
\ar@{->}[r]
\; & \; {\bf G} \ltimes_\red \ind_{\bf H}^{\bf G} \res^{\bf G}_{\bf H}(B) 
}
$$
in $ KK^{\check{\bf G}} $, where the unlabeled arrow is the canonical inclusion. 
The vertical comparison maps are given by the identity, and the covariant Morita 
equivalence $ {\bf G} \ltimes_\red \ind_{\bf H}^{\bf G} \res^{\bf G}_{\bf H}(B) \sim_M {\bf H} \ltimes_\red \res^{\bf G}_{\bf H}(B) $, respectively. 
Since the latter composition is the canonical embedding 
$ {\bf H} \ltimes_\red \res^{\bf G}_{\bf H}(B) \rightarrow {\bf G} \ltimes_\red \ind_{\bf H}^{\bf G} \res^{\bf G}_{\bf H}(B) $ 
we conclude that $ \res(\kappa_B) \circ \eta_{\res(B)} $ is the identity in $ KK^{\bf H} $. \qed

\section{The categorical Baum-Connes assembly map} \label{secbc} 

In this section we define and study the categorical assembly map for complex quantum groups in the framework of Meyer and Nest. 
Throughout we fix a positive real number $ q \neq 1 $ and let $ {\bf G} $ be the continuous field of locally compact quantum 
groups over $ [0,1] $ with fibers $ {\bf G}_t = G_{q^t} $ as described at the end of section \ref{seccflcqg}. 

From the construction of $ {\bf G} $ one sees that the maximal compact quantum subgroups in each fiber define a continuous field of 
compact quantum groups $ {\bf K} \subset {\bf G} $, and in fact a continuous field of compact open quantum subgroups. 
In particular, we have associated induction and restriction functors.  
Let us define full subcategories $ \CC, \CI \subset KK^{\bf G} $ by 
$$
\CC = \{A \in KK^{\bf G} \mid \res^{\bf G}_{\bf K}(A) \cong 0 \text{ in } KK^{\bf K} \} 
$$
and 
$$
\CI = \{A \in KK^{\bf G} \mid A \cong \ind^{\bf G}_{\bf K}(B) \text{ for some } B \in KK^{\bf K} \}.  
$$
Then the category $ \CC $ is localising, and we let $ \bra \CI \ket \subset KK^{\bf G} $ be the localising subcategory generated by $ \CI $. 

Recall that a pair $ (\mathfrak{L}, \mathfrak{N}) $ of thick subcategories 
of a triangulated category $ \T $ is called complementary if $ \T(L, N) = 0 $ 
for all $ L \in \mathfrak{L} $ and $ N \in \mathfrak{N} $, and every object $ A $ of $ \T $ fits into an exact triangle 
$ L \rightarrow A \rightarrow N \rightarrow L[1] $ 
with $ L \in \mathfrak{L} $ and $ N \in \mathfrak{N} $. 

\begin{prop} \label{complementarity}
The pair of categories $ (\bra \CI \ket, \CC) $ is complementary in $ KK^{\bf G} $. 
\end{prop} 

\proof Using the fact that $ {\bf K} \subset {\bf G} $ is open this is proved in the same way as Theorem 7.3 in \cite{Meyerhomalg2}. 
More precisely, since $ \ind_{\bf K}^{\bf G} $ is left adjoint to the restriction functor $ \res_{\bf K}^{\bf G} $ by Theorem \ref{opensubgroupinduction}, 
one finds that the ideal $ \ker(\res^{\bf G}_{\bf K}) $ in $ KK^{\bf G} $ has enough projective objects, and that these projective objects are 
precisely the direct summands of objects in $ \CI $.
According to Theorem 3.31 and Theorem 3.21 in \cite{Meyerhomalg2} this yields the claim. \qed   

Due to Proposition \ref{complementarity} each object $ A \in KK^{\bf G} $ fits into an exact 
triangle $ L_A \rightarrow A \rightarrow N_A \rightarrow L_A[1] $, and complementarity ensures that this triangle is unique up to isomorphism. 
This way one obtains a functor $ L: KK^{\bf G} \rightarrow KK^{\bf G} $ which is given by $ L(A) = L_A $ on objects, together 
with a natural transformation $ L(A) \rightarrow A $. The resulting morphism will also be referred to as a $ \CI $-cellular approximation, 
or simply cellular approximation, of $ A $. 

By definition, the assembly map for $ {\bf G} $ with coefficients in $ A \in KK^{\bf G} $ is the induced map 
$$
\mu^A: K_*({\bf G} \ltimes_\red L(A)) \rightarrow K_*({\bf G} \ltimes_\red A) 
$$
on the $ K $-theory of the reduced crossed products. 

We are mainly interested in the assembly maps for the fiber quantum groups, which are constructed in exactly the same way. More precisely, 
we obtain subcategories $ \CC_{q^t}, \CI_{q^t} $ and $ \bra \CI_{q^t} \ket $ in $ KK^{{\bf G}_t} = KK^{G_{q^t}} $ by replacing $ {\bf G} $ 
and $ {\bf K} $ in the above constructions with the fiber quantum groups $ {\bf G}_t = G_{q^t} $ and $ {\bf K}_t = K_{q^t} $, respectively. 
The same proof as for Proposition \ref{complementarity} shows that the pair $ (\bra \CI_{q^t} \ket, \CC_{q^t}) $ in $ KK^{G_{q^t}} $ for $ t \in [0,1] $  
is complementary, and we obtain cellular approximation functors $ L_{q^t}: KK^{G_{q^t}} \rightarrow KK^{G_{q^t}} $. 

\begin{definition} 
Let $ G_q $ be a complex semisimple quantum group and let $ A $ be a $ G_q $-$ C^* $-algebra. The categorical Baum-Connes assembly for $ G_q $ with coefficients 
in $ A $ is the map 
$$
\mu^A_q: K_*(G_q \ltimes_\red L_q(A)) \rightarrow K_*(G_q \ltimes_\red A) 
$$
induced from a cellular approximation $ L_q(A) $ of $ A $ in $ KK^{G_q} $. 
\end{definition} 

We shall now construct a cellular approximation $ \P_q \rightarrow \mathbb{C} $ for the trivial action of $ G_q $ starting from the Koszul complex 
of the representation ring $ R(K_q) = R(K) $. 
For later use, it will in fact be convenient to consider this construction on the level of the continuous field $ {\bf G} $. 

As a first step let us give a concrete description of the homogeneous space $ \ind_{\bf K}^{\bf G}(C[0,1]) = C_0({\bf G}/{\bf K}) $.  
We can identify 
$$ 
C_0({\bf G}/{\bf K}) = {\bf W}^* (1 \otimes C^*({\bf K})) {\bf W} \subset C({\bf K}) \otimes_{C[0,1]} C^*({\bf K}) = C_0({\bf G}).  
$$ 
In other words, $ C_0({\bf G}/{\bf K}) $ is isomorphic to $ C^*({\bf K}) $ in such a way that the left coaction of $ C^*({\bf K}) $ on $ C_0({\bf G}/{\bf K}) $ 
becomes the regular coaction. The left coaction of $ C({\bf K}) $ on $ C_0({\bf G}/{\bf K}) $ is identified 
with the fibrewise conjugation action of $ {\bf K} $ on $ C^*({\bf K}) $, implemented by $ {\bf W} \in M(C({\bf K}) \otimes_{C[0,1]} C^*({\bf K})) $. 

\begin{lemma} \label{kkgqhomogeneous}
We have 
$$
KK^{\bf G}(C_0({\bf G}/{\bf K}), C[0,1]) \cong R(K).  
$$
\end{lemma} 

\proof According to Theorem \ref{opensubgroupinduction} we obtain isomorphisms
\begin{align*}
KK^{\bf G}_*(C_0({\bf G}/{\bf K}), C[0,1]) &\cong KK^{\bf K}_*(C[0,1], C[0,1]) \cong K_*(C^*({\bf K})) = R(K), 
\end{align*}
using that the field $ C^*({\bf K}) $ is constant. \qed 

Recall that the representation ring $ R(K) = \mathbb{Z}[\varpi_1, \dots, \varpi_N] $ is a polynomial ring generated by the classes of the fundamental 
representations $ V(\varpi_1), \dots, V(\varpi_N) $ of $ K $. 
For each $ 1 \leq j \leq N $ we obtain a chain complex of free $ \mathbb{Z}[\varpi_j] $-modules 
\begin{center}
\begin{tikzcd}
      C^{(j)}: 0
			\arrow[r, ""]
			&
      \mathbb{Z}[\varpi_j]			
			\arrow[r, "\varpi_j - d_j"]
			&
			\mathbb{Z}[\varpi_j]
  		\arrow[r, ""]
			&
      0
\end{tikzcd}
\end{center}
Here the differential is given by multiplication with $ \varpi_j - d_j $, where $ d_j = \dim(V(\varpi_j)) $.
We form the Koszul complex $ C $ by taking the tensor product of these complexes, that is, 
$$ 
C = C^{(1)} \otimes_\mathbb{Z} \cdots \otimes_\mathbb{Z} C^{(N)}, 
$$
compare section 4.5 in \cite{Weibelhomalg}. By construction, the Koszul complex is a complex of free modules 
over $ R(K) = \mathbb{Z}[\varpi_1, \dots, \varpi_N] $.  
Writing $ e_i $ for the generator of the module $ \mathbb{Z}[\varpi_j] $ in degree $ 1 $ of $ C^{(j)} $, 
we may view the underlying module of $ C $ as the exterior algebra 
of the free $ R(K) $-module $ \Omega = \Omega^1_{\mathbb{Z}}(R(K)) $ 
with basis $ e_1, \dots, e_N $. More precisely, we get a $ R(K) $-basis of the component $ C_k $ of $ C $ in degree $ k \geq 1 $ given by the elements 
$$
e_{i_1} \wedge \cdots \wedge e_{i_k} 
$$
for $ 1 \leq i_1 < \dots < i_k \leq N $, such that the differential in $ C $ reads 
$$
\partial(e_{i_1} \wedge \cdots \wedge e_{i_k}) 
= \sum_{j = 1}^k (-1)^{j - 1} (\varpi_{i_j} - d_{i_j}) e_{i_1} \wedge \cdots \wedge e_{i_{j - 1}} \wedge e_{i_{j + 1}} \wedge \cdots \wedge e_{i_k}.  
$$
Since the elements $ \varpi_j - d_j $ form a regular sequence in $ \mathbb{Z}[\varpi_1, \dots, \varpi_N] $, the Koszul complex can be augmented to an exact 
chain complex 
\begin{center}
\begin{tikzcd}
			0
			\arrow[r, ""]
			&
			\Lambda^N \Omega  
			\arrow[r, "\partial"]
			&
			\Lambda^{N - 1} \Omega
			\arrow[r, "\partial"]
			&
			\cdots 
			\arrow[r, "\partial"]
			&
			\Lambda^1 \Omega 
  		\arrow[r, "\partial"]
			&
			\Lambda^0 \Omega 
  		\arrow[r, "\epsilon"]
			&
			\mathbb{Z} 
			\arrow[r, ""]
			&
			0
\end{tikzcd}
\end{center}
of $ R(K) $-modules. Here $ \epsilon: \Lambda^0 \Omega = R(K) \rightarrow \mathbb{Z} $ is the augmentation homomorphism given by $ \epsilon(\varpi_j) = d_j $. 
Note that the resulting augmented chain complex is free as a complex of $ \mathbb{Z} $-modules, so that we can find a $ \mathbb{Z} $-linear contracting 
homotopy for it. 

We shall now lift this complex to the level of $ KK^{\bf G} $. For $ 0 \leq k \leq N $ let $ r_k = \binom{N}{k} $ be the rank 
of the $ R(K) $-module $ \Lambda^k \Omega $ and consider the $ {\bf G} $-$ C^* $-algebra 
$$ 
P_k = \bigoplus_{j = 1}^{r_k} C_0({\bf G}/{\bf K}),  
$$ 
which we shall view as a direct sum of copies of $ C_0({\bf G}/{\bf K}) $ indexed by the basis vectors $ e_{i_1} \wedge \cdots \wedge e_{i_k} $ from above. 
Then we may identify $ K_0(P_k) = \Lambda^k \Omega $, and we claim that we can build the differential $ \partial $ 
of the Koszul complex using morphisms in $ KK^{\bf G} $. To this end it is enough to implement the multiplication action of $ \varpi_j $ on $ R(K) $ 
by an element of $ KK^{\bf G}(C_0({\bf G}/{\bf K}), C_0({\bf G}/{\bf K})) $. 
The comultiplication $ \hat{\Delta}: C^*({\bf K}) \rightarrow M(C^*({\bf K}) \otimes_{C[0,1]} C^*({\bf K})) $, 
composed with the representation $ C^*({\bf K}) \rightarrow \KH(V(\varpi_j)) $ provides a representative of the desired class, 
using the $ {\bf G} $-equivariant Morita equivalence between $ C^*({\bf K}) \otimes \KH(V(\varpi_j)) $ and $ C^*({\bf K}) = C_0({\bf G}/{\bf K}) $. 

The map $ \epsilon: R(K) \rightarrow \mathbb{Z} $ can be lifted to the class induced by 
$$ 
{\bf G} \ltimes_\red C_0({\bf G}/{\bf K}) \sim_M C^*_\red({\bf K}) \rightarrow C^*_\red({\bf G}) \cong {\bf G} \ltimes_\red C[0,1]
$$
under Baaj-Skandalis duality. Using Lemma \ref{kkgqhomogeneous} one checks that this yields a complex 
\begin{center}
\begin{tikzcd}
			0
			\arrow[r, ""]
			&
			P_N
			\arrow[r, "\partial"]
			&
			P_{N - 1} 
			\arrow[r, "\partial"]
			&
			\cdots 
			\arrow[r, "\partial"]
			&
			P_1 
  		\arrow[r, "\partial"]
			&
			P_0
  		\arrow[r, "\epsilon"]
			&
			C[0,1]
			\arrow[r, ""]
			&
			0
\end{tikzcd}
\end{center}
in $ KK^{\bf G} $, which upon taking $ K $-theory gives the augmented Koszul complex. 
Here, by slight abuse of notation, we denote the boundary maps at the level of $ KK^{\bf G} $ by the same symbols as their 
counterparts in the Koszul complex. 

\begin{lemma} \label{resolution}
The complex $ P $ defines a $ \CI $-projective resolution of the trivial $ {\bf G} $-$ C^* $-algebra $ C[0,1] $ in $ KK^{\bf G} $. 
\end{lemma} 

\proof By construction the objects $ P_n $ are $ \CI $-projective. For exactness it suffices to observe
that the action of $ {\bf K} $ on $ C_0({\bf G}/{\bf K}) $ is $ K $-theoretically trivial, 
so that we can write down a contracting homotopy at the level of $ KK^{\bf K} $, by lifting a $ \mathbb{Z} $-linear 
contracting homotopy of the augmented Koszul complex. \qed 

From Lemma \ref{resolution} we obtain a $ \CI $-cellular approximation of $ C[0,1] $ using general machinery from homological algebra in triangulated 
categories \cite{Meyerhomalg2}. More precisely, starting from $ P $ one constructs, in the terminology of \cite{Meyerhomalg2}, a phantom castle 
over $ C[0,1] $ in $ KK^{\bf G} $. The corresponding cellular approximation tower yields an exact 
triangle $ \P \rightarrow C[0,1] \rightarrow N \rightarrow \P[1] $ in $ KK^{\bf G} $
such that $ \P \in \bra \CI \ket $ and $ N \in \CC $. The object $ \P $ can be described as homotopy colimit of the cellular approximation tower, 
and together with the morphism $ \P \rightarrow C[0,1] $ in $ KK^{\bf G} $ from the triangle this provides a $ \CI $-cellular approximation of $ C[0,1] $. 

The analogous constructions for the fiber quantum groups provide cellular approximations as well. 
More precisely, for every $ t \in [0,1] $, the morphism $ \P_{q^t} \rightarrow \mathbb{C} $ in $ KK^{G_{q^t}} $ 
induced from $ \P \rightarrow C[0,1] $ in $ KK^{\bf G} $ is a $ \CI_{q^t} $-cellular approximation of $ \mathbb{C} \in KK^{G_{q^t}} $. 

We are mainly interested in the case $ t = 1 $ and the resulting assembly map for $ G_q $, but we shall first 
consider the case $ t = 0 $, corresponding to $ G_1 = K \bowtie \hat{K} $. For this purpose we need an equivariant version of Baaj-Skandalis duality
for the group $ K $. Let $ K \ltimes K $ be the compact group obtained as the semidirect product of $ K $ acting on itself by conjugation. 

\begin{prop} \label{relBSduality}
Taking reduced crossed products with respect to $ K $ and $ \hat{K} $, respectively, induces equivalences of triangulated categories
\begin{align*} 
J_K^K&: KK^{K \ltimes K} \rightarrow KK^{K \bowtie \hat{K}} \\
J_{\hat{K}}^K&: KK^{K \bowtie \hat{K}} \rightarrow KK^{K \ltimes K} 
\end{align*} 
which are mutually inverse up to natural isomorphism. 
\end{prop} 

\proof Let $ A \in KK^{K \ltimes K} $ and denote by $ \alpha: A \rightarrow C(K) \otimes A $ the coaction corresponding to the action of 
the second factor in the semidirect product. We equip the crossed product $ K \ltimes_\red A = [(C^*(K) \otimes 1)\alpha(A)] \subset \LH(L^2(K) \otimes A) $ 
with the $ K $-action given by conjugation on $ C^*(K) \subset M(K \ltimes_\red A) $ and the action of the first factor of $ K \ltimes K $ 
on $ A \subset M(K \ltimes_\red A) $. 
This is implemented by the tensor product of the conjugation action of $ K $ on $ L^2(K) $ and the action of the first factor of $ K \ltimes K $ on $ A $ 
inside $ L^2(K) \otimes A $. 
Together with the dual action of $ \hat{K} $ this yields a $ K \bowtie \hat{K} $-$ C^* $-algebra structure on $ K \ltimes_\red A = J_K^K(A) $, 
compare Proposition 3.2 in \cite{NVpoincare}. 
Using the universal property of equivariant $ KK $-theory, see Theorem \ref{kkuniversal}, one verifies that this construction extends to a triangulated 
functor $ J_K^K: KK^{K \ltimes K} \rightarrow KK^{K \bowtie \hat{K}} $ as stated. 

If $ B \in KK^{K \bowtie \hat{K}} $ and $ \beta: B \rightarrow M(C^*(K) \otimes B) $ denotes the coaction corresponding to the 
action of $ \hat{K} \subset K \bowtie \hat{K} $, then $ \hat{K} \ltimes B = [(C(K) \otimes 1) \beta(B)] \subset \LH(L^2(K) \otimes B) $ 
is naturally equipped with the dual action of $ K $. We obtain an additional action of $ K $ by conjugation with the tensor product action on the Hilbert 
module $ L^2(K) \otimes B $ given by conjugation in the first factor and the given $ K $-action on $ B $ in the second. 
In the same way as above one checks that these actions combine to an action of $ K \ltimes K $ on $ \hat{K} \ltimes B = J_{\hat{K}}^K(B) $, 
and that this construction induces a triangulated functor $ J_{\hat{K}}^K: KK^{K \bowtie \hat{K}} \rightarrow KK^{K \ltimes K} $. 

In order to show that the functors $ J_K^K $ and $ J_{\hat{K}}^K $ are mutually inverse up to natural isomorphism it suffices to inspect the proof of 
Takesaki-Takai duality, see \cite{BSUM} or chapter 9 in \cite{Timmermannbook}, noting that the multiplicative unitary $ W \in \LH(L^2(K) \otimes L^2(K)) $ 
for the group $ K $ commutes with the diagonal action of $ K $ on $ L^2(K) \otimes L^2(K) $ induced by conjugation in both factors. \qed 

The group $ C^* $-algebra $ C^*_\red(K \bowtie \hat{K}) $ of $ G_1 = K \bowtie \hat{K} $ can be identified with 
$$
C^*_\red(K \bowtie \hat{K}) = [(C^*(K) \otimes 1)\ad(C(K))] = K \ltimes_\ad C(K),
$$
the crossed product of $ C(K) $ with respect to the adjoint action of $ K $. Explicitly, the corresponding coaction 
$ \ad: C(K) \rightarrow C(K) \otimes C(K) $ is given by $ \ad(f)(s,t) = f(sts^{-1}) $. 
Note that $ \ad(f) = U^*(1 \otimes f) U $ where $ U \in \LH(L^2(K) \otimes L^2(K)) $ is the unitary operator given by 
$$
U(h \otimes k)(s,t) = h(s) k(s^{-1} t s).  
$$
With this description at hand let us explain how to interpret reduced crossed products with respect to $ K \bowtie \hat{K} $ as iterated crossed products, 
first by $ \hat{K} $ and then by $ K $. More precisely, let $ A \in KK^{K \bowtie \hat{K}} $. Then the 
coaction $ \gamma: A \rightarrow M(C_0(K \bowtie \hat{K}) \otimes A) $ corresponds to a pair of 
coactions $ \alpha: A \rightarrow C(K) \otimes A, \lambda: A \rightarrow M(C^*(K) \otimes A) $ 
such that $ \gamma(a) = (\id \otimes \lambda) \alpha $, see again Proposition 3.2 in \cite{NVpoincare}. We therefore obtain 
\begin{align*}
(K \bowtie \hat{K}) &\ltimes_\red A = [(C^*(K \bowtie \hat{K}) \otimes 1) \gamma(A)] \\
&= [(C^*(K) \otimes 1 \otimes 1) U_{12}^* (1 \otimes C(K) \otimes 1) U_{12} (\id \otimes \lambda)\alpha(A)] \\
&= [(C^*(K) \otimes 1 \otimes 1) U_{12}^* (1 \otimes C(K) \otimes 1) U_{12} W^*_{12} (\sigma \otimes \id)(\id \otimes \alpha)\lambda(A) W_{12}] \\
&= [(C^*(K) \otimes 1 \otimes 1) U_{12}^* (1 \otimes C(K) \otimes 1) (\sigma \otimes \id)(\id \otimes \alpha)\lambda(A) U_{12}] \\
&\cong K \ltimes_{\red, \tilde{\alpha}} \hat{K} \ltimes_\red A,
\end{align*} 
where $ \tilde{\alpha}: \hat{K} \ltimes_\red A \rightarrow C(K) \otimes (\hat{K} \ltimes_\red A) $ is the action given by 
\begin{align*} 
\tilde{\alpha}((f \otimes 1) \lambda(a)) &= U_{12}^* (\sigma \otimes \id)(\id \otimes \alpha)((f \otimes 1)\lambda(a)) U_{12} \\
&= (\ad(f) \otimes 1)(\id \otimes \lambda)\alpha(a). 
\end{align*}

With these preparations in place we are ready to study the assembly map for the quantum group $ G_1 $. 

\begin{theorem} \label{BCG1} 
The Drinfeld double $ G_1 = K \bowtie \hat{K} $ of the classical group $ K $ satisfies the strong Baum-Connes property, that is, we 
have $ KK^{G_1} = \bra \CI_1 \ket $. 
\end{theorem} 

\proof For $ A \in KK^K $ we denote by $ \res^K_{K \ltimes K}(A) \in KK^{K \ltimes K} $ the $ K \ltimes K $-$ C^* $-algebra obtained by letting $ K \ltimes K $ 
act via the group homomorphism $ K \ltimes K \rightarrow K $ given by projection to the first factor. 
Using the fact that we may write crossed products with respect to $ K \bowtie \hat{K} $ as iterated crossed products as explained above, we obtain equivariant 
isomorphisms and Morita equivalences 
\begin{align*}
(K \bowtie \hat{K}) \ltimes_\red J_K^K(\res^K_{K \ltimes K}(A)) &\cong (K \bowtie \hat{K}) \ltimes_\red (C^*(K) \otimes A) \\
&\cong K \ltimes_\red (\KH(L^2(K)) \otimes A) \\
&\sim_M K \ltimes_\red A \\
&\sim_M (K \bowtie \hat{K}) \ltimes_\red \ind_K^{K \bowtie \hat{K}}(A) 
\end{align*}
with respect to the action of the dual of $ K \bowtie \hat{K} $. This yields an isomorphism 
$$ 
J_K^K(\res^K_{K \ltimes K}(A)) \cong \ind_K^{K \bowtie \hat{K}}(A) 
$$ 
in $ KK^{K \bowtie \hat{K}} $. 

Taking into account Proposition \ref{relBSduality}, this observation shows that it is enough to 
verify $ KK^{K \ltimes K} = \bra \res_{K \ltimes K}^K(A) \mid A \in KK^K \ket $. This, in turn, is an immediate consequence of the strong Baum-Connes-property for 
the dual of the compact Lie group $ K \ltimes K $, see \cite{MNcompact}. 
In fact, it suffices to consider $ C^* $-algebras $ A $ with trivial $ K $-action on the right hand side. \qed 

Let us now formulate and prove our main result. 

\begin{theorem} \label{categoricalBCcqg}
The complex quantum group $ G_q $ satisfies the Baum-Connes conjecture. That is, the assembly map 
$$
\mu_q: K_*(G_q \ltimes_\red \P_q) \rightarrow K_*(C^*_\red(G_q)) 
$$
is an isomorphism. 
\end{theorem} 

\proof The map $ \mu_q $ fits into the commutative diagram 
\begin{center}
\begin{tikzcd}
      K_*(G_q \ltimes_\red \P_q) 
			\arrow[r, "\mu_q"]
			&
			K_*(C^*_\red(G_q)) \\
      K_*({\bf G} \ltimes_\red \P) 
			\arrow[r, "\mu"]
			\arrow[d, ""]
			\arrow[u, ""]
			&
			K_*(C^*_\red({\bf G})) 
			\arrow[d, ""]
			\arrow[u, ""] \\
 	    K_*(G_1 \ltimes_\red \P_1) 		
		  \arrow[r, "\mu_1"]
			&
			K_*(C^*_\red(G_1)) 
\end{tikzcd}
\end{center}
obtained by specialising the assembly map for $ {\bf G} $ to the fibers at the endpoints of the interval $ [0,1] $. 
According to Theorem 6.2 and Theorem 6.5 in \cite{MVcqgdeformation}, the right vertical maps in this diagram are both isomorphisms. 

Consider next the left hand side of the diagram. Let $ F $ be the homological functor from $ KK^{\bf G} $ to the category of abelian groups given 
by $ F(A) = K({\bf G} \ltimes_\red A) $. We may compute $ K_n({\bf G} \ltimes_\red \P) $ using the ABC spectral sequence from \cite{Meyerhomalg2} for the 
functor $ F $. The first page of this spectral sequence can be written as $ E^1_{mn} = F_n(P_m) $, where $ P $ is the projective resolution in $ KK^{\bf G} $ 
obtained further above in this section and $ F_n(A) = K_n({\bf G} \ltimes_\red A) $. 
According to Proposition 4.10 and Proposition 3.28 in \cite{Meyerhomalg2}, the ABC spectral sequence collapses after $ N + 1 $ steps and converges 
towards the $ K $-theory of $ {\bf G} \ltimes_\red \P $. 

In the same way we obtain homological functors 
and their associated ABC spectral sequences with respect to the fibers of the continuous field $ {\bf G} $. 
Since the continuous fields $ {\bf G} \ltimes_\red P_m $ in the complex $ {\bf G} \ltimes_\red P $ are Morita equivalent to constant fields, 
evaluation at the endpoints of $ [0,1] $ induces isomorphisms of the $ E^1 $-terms of these spectral sequences. 
It follows that the left vertical maps in the above diagram are both isomorphisms. 

The strong Baum-Connes property for $ G_1 $ obtained in Theorem \ref{BCG1} implies that the bottom horizontal map is an isomorphism as well. 
Combining these facts we deduce that $ \mu $ and $ \mu_q $ are isomorphisms. \qed 

From the proof of Theorem \ref{categoricalBCcqg} it follows in particular that the assembly map 
constructed using deformations in \cite{MVcqgdeformation} is canonically isomorphic to the assembly map $ \mu_q $ 
obtained from the categorical setup discussed here. 

Let us remark that braided tensor products allow one to describe the assembly map for $ G_q $ with arbitrary coefficients using the cellular 
approximation $ \P_q $. More precisely, Theorem 3.6 in \cite{NVpoincare} with $ G = G_q $ and $ H = K_q $ yields 
$$ 
C_0(G_q/K_q) \boxtimes_{G_q} A = \ind_{K_q}^{G_q}(\mathbb{C}) \boxtimes_{G_q} A \cong \ind_{K_q}^{G_q}(\mathbb{C} \boxtimes_{K_q} \res^{G_q}_{K_q}(A)) 
= \ind_{K_q}^{G_q} \res^{G_q}_{K_q}(A)
$$ 
for all $ A \in KK^{G_q} $. It follows that $ \P_q \boxtimes_{G_q} A $ is contained in $ \bra \CI_q \ket $, and the canonical 
morphism $ \P_q \boxtimes_{G_q} A \rightarrow \mathbb{C} \boxtimes_{G_q} A \cong A $ is a cellular approximation of $ A $. 
In particular, the Baum-Connes assembly map for $ G_q $ with coefficients in $ A \in KK^{G_q} $ is the 
induced map $ K_*(G_q \ltimes_\red (\P_q \boxtimes_{G_q} A)) \rightarrow K_*(G_q \ltimes_\red A) $ in $ K $-theory. 
According to Theorem \ref{BCG1}, this map is always an isomorphism for $ q = 1 $.

\bibliographystyle{hacm}

\bibliography{cvoigt}

\end{document}